\documentclass[12pt,a4paper]{article}

\usepackage{amsmath}

\title{Asymptotic behavior and halting probability of Turing 
Machines\footnote{First version: March 2006. Article in press. 
{\it Chaos, Solitons \& Fractals}~(2006), doi:10.1016/j.chaos.2006.08.022.}}

\author{Germano D'Abramo\vspace{0.3cm}\\
{\small Istituto Nazionale di Astrofisica,}\\
{\small Via Fosso del Cavaliere n.~100,}\\
{\small 00133 Roma, Italy.}\\
{\small E--mail: {\tt Germano.Dabramo@iasf-roma.inaf.it}}}

\begin{document}

 \date{}

\maketitle

\begin{abstract}

Through a straightforward Bayesian approach we show that under some general 
conditions a maximum running time, namely the number of discrete steps 
performed by a computer program during its execution, can be defined such 
that the probability that such a program will halt {\em after} that time is 
smaller than any arbitrary fixed value. Consistency with known results and 
consequences are also discussed.

\end{abstract}

\section{Introductory remarks}

As it has been proved by Turing in 1936~\cite{turi}, if we have a 
program $p$ running on an Universal Turing Machine (UTM), then we have no 
general, finite and deterministic algorithm which allows us to know whether 
and when it will halt (this is the well known {\em halting problem}). 
That is to say that the halting behavior of a program, with the trivial 
exception of the simplest ones, is not computable and predictable by a 
unique, general procedure. 

In this paper we show that, for what concerns the probability of its halt, 
every program running on an UTM is characterized by a peculiar asymptotic 
behavior in time. Similar results have been obtained by 
Calude~{\em et al.}~\cite{calude1} and by Adamyan~{\em et al.}~\cite{ada} 
through a different approach, which makes use of quantum computation.

\section{Probabilistic approach}

Given a program $p$ of $n$~bits, it is always possible to slightly change 
its code (increasing its size by a small, fixed amount of bits, let's say 
of $s$~bits with $s\ll n$) in order to include a 
progressive, integer counter that starts counting soon after $p$ starts to 
run on an UTM and stops, printing the total number of steps done by the 
UTM during the execution of $p$, soon after $p$ halts. Let us call this 
new program $p'$. Its algorithmic size is then of $n+s$~bits.

As it has been said in the previous Section, for an arbitrary program of 
$n$~bits no general, finite and deterministic procedure exists that 
allows us to know whether such program will ever halt or will keep running 
forever on an UTM~\cite{turi}. Thus, in our case, we have no finite, 
deterministic algorithmic  procedure to decide whether and when $p'$ 
will print the integer number of steps done by the UTM during the execution 
of $p$, until the halt.

Let us now make a brief digression on algorithmic complexity. Suppose we 
have a (randomly chosen) binary string of $n$~bits. Its algorithmic 
complexity (for a rigorous definition of algorithmic complexity, 
see~\cite{chai1}) is less than or equal to $n+c$~bits, where $c$ is a 
specified constant in the chosen language. The case of a self-delimiting UTM 
(see~\cite{chai1}), for which the algorithmic complexity of a program of 
$n$~bits is less than or equal to $n+O(\log_2 n)$~bits, is dealt with later.
The {\em a priori} probability $P_1$ that its algorithmic complexity  is 
equal to $k$~bits\footnote{And thus, following the definition of algorithmic 
complexity~\cite{chai1}, the {\em a priori} probability that a minimum program 
of $k$~bits exists that generates the string of $n$~bits as output.}, with 
$k\leq n+c$, is:

\begin{equation} 
P_1=\frac{2^k}{2^{n+c+1}-2},
\label{eq1}
\end{equation}
since the total number of possible generating strings less than or equal to 
$n+c$~bits in size is $2^{n+c+1}-2$, while there are only $2^k$ strings of 
$k$~bits.

Now we want to calculate the somewhat different {\em a priori} probability 
$P_2$ that a (randomly chosen) program of algorithmic complexity of $k$~bits 
that produces an output\footnote{Note that this condition is crucial.} 
generates just a string of $n$~bits, with $n+c\geq k$.

From Bayes theorem we have that such probability can be calculated as the 
ratio between $P_1$ of eq.~(\ref{eq1}) and the sum of all the probability 
$P_1$, with $n+c$ ranging from $k$ to infinity:

\begin{equation} 
P_2=\frac{\frac{2^k}{2^{n+c+1}-2}}{\sum_{i=k}^\infty 
\frac{2^k}{2^{i+1}-2}}.
\label{eq2}
\end{equation}

The probability $P_2$ can be easily rewritten and simplified as follows:

\begin{equation} 
P_2=\frac{1}{2^{n+c+1}-2}\cdot\frac{1}{\sum_{i=k}^\infty 
\frac{1}{2^{i+1}-2}}\simeq\frac{2^k}{2^{n+c+1}-2}\simeq\frac{1}{2^{n+c-k+1}}
\,\,\,\,\,\,\, n+c\geq k,
\label{eq3}
\end{equation}
since the sum $\sum_{i=k}^\infty \frac{1}{2^{i+1}-2}$ is safely approximable
to $1/2^k$ for $k\geq 10$.

Suppose now that we have a program $p'$ of algorithmic complexity of $k$~bits 
that will halt for sure, even if we are not able to know if and when it will 
do so by simply inspecting its code, and that contains a discrete step counter 
which counts the total number of steps done by the program until the halt and 
prints it as output\footnote{For every size $k$ we can add few instructions to 
implement the counter, altering only slightly its length, as explained at the
beginning of this Section.}. 
Applying Bayes theorem again, the probability that the output string is of a 
size greater than or equal to $m$~bits (with $m+c\geq k$) is:

\begin{equation} 
P(\textrm{\small size}\geq m)=\frac{\sum_{s=m}^\infty
\frac{2^k}{2^{s+c+1}-2}}{\sum_{i=k}^\infty 
\frac{2^k}{2^{i+1}-2}},\,\,\,\,m+c\geq k.
\label{eq4}
\end{equation}

Using the same approximation done in eq.~(\ref{eq3}) the above equation 
becomes:

\begin{equation} 
P(\textrm{\small size}\geq m)\simeq\frac{1}{2^{m+c-k}},\,\,\,\,m+c\geq k.
\label{eq4b}
\end{equation}

Thus, the probability that the output string is less than $m$~bits in size
is:

\begin{equation} 
P(\textrm{\small size}<m)=1-\frac{\sum_{s=m}^\infty
\frac{2^k}{2^{s+c+1}-2}}{\sum_{i=k}^\infty 
\frac{2^k}{2^{i+1}-2}}\simeq1-\frac{1}{2^{m+c-k}},\,\,\,\,m+c\geq k.
\label{eq5}
\end{equation}

It is easy to verify through a numerical check that for $m+c\geq k + 50$, 
and thus even better for $m\geq k + 50$, the  value of 
$P(\textrm{\small size}<m)$ is practically equal to $1$.

This means that {\em if} a program $p$ of algorithmic complexity of 
$k$~bits halts, then it will do so in such a way that the total number 
of steps is a number of  size less than or equal to $k+50$~bits with very 
high probability. Thus, the decimal number of steps $t$ done by the program 
before the halt (the discrete time before the halt, or the {\em 
characteristic time}) can not be more than $2^{k+51}$, with very high 
probability.

The above analysis and results hold also if the program $p$ has a size of
$k$~bits, but a smaller algorithmic complexity. In such cases the 
characteristic time $2^{k+51}$ provides a more conservative halting time.  

If we consider a self-delimiting UTM, then the algorithmic complexity of a 
program of $n$~bits is less than or equal to $n+O(\log_2 n)$~bits 
(see~\cite{chai1}). With such complexity measure eq.~(\ref{eq1}) becomes

\begin{equation} 
P_1=\frac{2^k}{2^{n+O(\log_2 n)+1}-2},
\label{eq1b}
\end{equation}
and eq.~(\ref{eq2}) consequently becomes 

\begin{equation} 
P_2=\frac{\frac{2^k}{2^{n+O(\log_2 n)+1}-2}}{\sum_{i=l}^\infty 
\frac{2^k}{2^{i+O(\log_2 i)+1}-2}},
\label{eq2b}
\end{equation}
where $l$ is an integer such that $l+O(\log_2 l)=k$.

Equation~(\ref{eq5}) then becomes

\begin{equation} 
P(\textrm{\small size}<m)=1-\frac{\sum_{s=m}^\infty
\frac{2^k}{2^{s+O(\log_2 s)+1}-2}}{\sum_{i=l}^\infty 
\frac{2^k}{2^{i+O(\log_2 i)+1}-2}},\,\,\,\,m+O(\log_2 m)\geq k.
\label{eq5b}
\end{equation}

If we consider cases with $l\gg O(\log_2 l)$ and $m\gg O(\log_2 m)$, then 
eq.(\ref{eq5b}) can be safely approximable by

\begin{equation} 
P(\textrm{\small size}<m)\simeq 1-\frac{\sum_{s=m}^\infty
\frac{2^k}{2^{s+1}-2}}{\sum_{i=k}^\infty 
\frac{2^k}{2^{i+1}-2}}\simeq 1-\frac{1}{2^{m-k}},\,\,\,\,m\geq k.
\label{eq5c}
\end{equation}
Thus, for $l\gg O(\log_2 l)$, where $l+O(\log_2 l)=k$, and $m\gg O(\log_2 m)$,
$P(\textrm{\small size}<m)\simeq 1-2^{k-m}$ and again the value of  
$P(\textrm{\small size}<m)$ is practically equal to $1$ for $m\geq k + 50$, as
before. 

Actually, it is possible to prove\footnote{Given eq.~(\ref{eq5}) or 
eq.~(\ref{eq5b}), it is always possible to find a constant $b>1$ such that:
$$ 1-\frac{\sum_{s=m}^\infty
\frac{2^k}{2^{s+1}-2}}{\sum_{i=\frac{k}{b}}^\infty 
\frac{2^k}{2^{bi+1}-2}}\simeq 1-\frac{1}{2^{m-b-k}}\leq  
P(\textrm{\small size}<m)\leq 1.$$
For $m-b\geq k + 50$, and thus for $m\geq k + 50 + b$, the  value of 
$P(\textrm{\small size}<m)$ is practically equal to $1$.} that for a 
suitable constant $b>1$ the function $1-2^{k-m+b}$ is a lower bound for the 
value of $P(\textrm{\small size}<m)$ for both complexity measures, and thus
for $m\geq k + 50 + b$, the  value of $P(\textrm{\small size}<m)$ becomes
closer and closer to $1$.

\section{Apparent paradoxes and their solution}

Now, consider the following peculiar case, which obviously provides a 
limitation in the application of eq.~(\ref{eq5}), as every peculiar case does
with probabilistic approaches. 

It is possible to write a simple program $p$ of $\lceil \log_2 n\rceil + 
c$~bits that, given the decimal number $n$,  recursively calculates the 
number $2^{n + 1}-2$, counts from $1$ to $2^{n + 1}-2$ and then stops. As it 
is obvious, the program $p$ runs for a total number of steps greater than 
$2^{n+ 1}-2$ and then stops, but eq.~(\ref{eq5}) may give an almost zero 
probability for such a result, since for suitable values of $n$ the 
characteristic time $2^{\lceil \log_2 n\rceil + c+51}$ can be much much 
smaller than $2^{n+ 1}-2$.

As it always happens with probabilistic treatments, if, analyzing the code, we 
are able to know in advance that the program will halt and we are also able 
to know when it will halt, then eq.~(\ref{eq5}) is of poor use. But this
does not dismiss our results as meaningless: it is obvious that when we have a
deterministic solution of a problem, the probabilistic one simply does not 
apply.

Another apparent paradox with the above results is the following. It would 
seem possible to write a relatively small program $h$, again of 
$\lceil \log_2 n\rceil + c$~bits and with $n\gg \lceil \log_2 n\rceil + c$, 
which lists all the $2^{n+1}-2$ programs of size less than or equal to 
$n$~bits, runs each program for a discrete time equal to $2^{k+51}$, where $k$ 
is the size of the program, and stores the output strings of the halting ones.
Then, simply printing as output a string greater by one unit than the greatest 
among those stored, the program $h$ is able to provide a string of 
algorithmic complexity greater than or equal to $n$~bits with a very high 
probability. This seems to challenge our results since 
$n\gg \lceil \log_2 n\rceil + c$ and $h$ would be able to print a string
of algorithmic complexity greater than its own size. But it does not. 

As a matter of fact, among all the program of size less than or equal to 
$n$~bits executed by $h$ there is the program $h$ itself, which we will call 
$h_2$, and this fact generates a contradiction that does not allow the 
program $h$ to produce any meaningful output. 

In fact, the size of $h_2$ is obviously equal to 
$\lceil \log_2 n\rceil + c$~bits and we already know that after a 
characteristic time of $2^{\lceil \log_2 n\rceil + c +51}$ steps it will 
be still running (the running time of $h$, and thus of $h_2$, is obviously 
greater than $2^{n+51}$ steps). But $h_2$ should halt by definition and 
thus we {\em surely} have among the programs of size less than or equal to $n$ 
selected by $h$ as non halting, an halting one, making the procedure 
$h$ meaningless. 

Besides, note that if we accept the output $s$ of $h$ as true, then the 
output of $h$ has to be at least $s+1$, since also the output of $h_2$ is 
equal to $s$, and so on, endlessly.

Stated in other words, it is not possible to use our results to write a 
mechanical procedure able to print a string more complex than the mechanical 
procedure itself, as it should be according to the definition of algorithmic 
complexity.

\section{Mathematical implications}

Being able to solve the halting problem has unimaginable mathematical 
consequences since many unanswered mathematical problems, such as the 
Goldbach's conjecture, the extended Riemann hypothesis and others, can be 
solved if one is able to show whether the program written to find a single 
finite counterexample will ever halt~\cite{chai2,gard}.

However, the estimate of the characteristic time done in the previous Section, 
namely $t\simeq 2^{n+51}$, where $n$ is the size of the program, 
shows that the result obtained in eq.~(\ref{eq5}) is not much 
useful for the practical resolution of the halting problem, even for a 
probabilistic one, since almost all the interesting programs have a size 
much greater than $50$~bits, giving astronomically huge characteristic 
times.

Anyway, our result should be of some theoretical interest since it shows 
an asymptotic behavior typical of every Turing Machine. All this seems to 
shed new light on the intrinsic significance of the algorithmic size or, 
more precisely, of the algorithmic complexity of a program encoding a 
mathematical problem. As a matter of fact, such a low-level and 
low-informative property of a program, as the number $n$ of its bit-size, 
seems to be strongly related to its halting behavior, and thus, according 
to the above-mentioned mathematical connections, it seems to be intimately 
linked to the high level, mathematical truth encoded in the program. 
Calude~{\em et al.}~\cite{calude2} have recently proposed a way to evaluate
the difficulty of a finitely refutable mathematical problem which is based 
just on the algorithmic complexity, in a fixed language, of the Turing 
Machine encoding the problem.

Consider the Riemann hypothesis, for instance~\cite{dusa}. If I am able to 
show that a program of $n$~bits, written to find a finite counterexample, 
will never halt with probability greater than $99.99999999\%$, then I may 
safely say that Riemann hypothesis is almost certainly true. The singular 
aspect here is that to be able to make such a claim I need only a finite 
number of numerical checks of the conjecture to reach a probability of 
$99.99999999\%$, out of an infinite number of zeros of the zeta function 
to be checked. Honestly speaking, it appears quite surprising. 

The above argument could also be seen as a re-proposition of the Humian 
induction problem~\cite{hume}, this time applied to finitely refutable 
mathematical statements: if a finitely refutable mathematical statement, 
encoded in a program of $n$~bits of which we are not able to know if it will 
halt simply inspecting its code, holds true for about $t\simeq 2^{n+51}$ 
steps, then it is definitely true with a fixed arbitrary high probability.
As a matter of fact, all this seems to give a strong quantitative support
to the inspiring principles of ``experimental mathematics'', proposed
with force by many scholars in the last years~\cite{chai3}.

\section*{Acknowledgments}

I wish to thank Prof.~Cristian Calude for helpful comments and suggestions.

\end{document}